\documentclass[a4paper]{article}
\newtheorem{theorem}{Theorem}[section]
\newtheorem{lemma}{Lemma}[section]

\newtheorem{proposition}{Proposition}[section]

\newenvironment{proof}{{\bf Proof.}}{\par\hspace{25em}\rule{1ex}{1ex}\vspace{1ex}}
\newenvironment{example}{{\bf Example.}}{\par\hspace{25em}\rule{1ex}{1ex}\vspace{1ex}}

\newcommand{\ds}{\displaystyle}

\title{On extrinsic geometry of unit normal vector fields of Riemannian
hyperfoliations.\footnote{Math. Publ. Debrecen. 2003,  63/4,   555-567}}
\author{Yampolsky A.}
\date{}
\begin{document}
\maketitle
\begin{abstract}
    We consider a unit normal vector field of (local) hyperfoliation on a given Riemannian
    manifold as a submanifold in the unit tangent bundle with Sasaki metric. We give an explicit
    expression of the second fundamental form for this submanifold and a rather simple condition
    its totally geodesic property in the case of a totally umbilic hyperfoliation.
    A corresponding example shows the non-triviality of this condition. In the 2-dimensional
    case, we give a complete description of Riemannian manifolds admitting a geodesic
    unit vector field with totally geodesic property.
    \\[1ex]
    {\it AMS subject class:} Primary 54C40,14E20; Secondary 46E25, 20C20
\end{abstract}

\section*{Introduction}

    Let $(M,g)$ be an $(n+1)$-dimensional Riemannian manifold  with metric $g$
    and $\xi$ a fixed unit vector field on $M$.
    Consider $\xi$ as a (local) mapping $\xi : M \to T_1M $. Then the image
    $\xi(M) $  is a submanifold in the unit tangent sphere bundle $T_1M$.
    The Sasaki metric on the tangent bundle $TM$ induces the Riemannian metric
    on $T_1M$ and on $\xi(M)$ as well.
    So, one may use notions from the geometry of submanifolds to determine
    geometrical characteristics of a unit vector field.

    A unit vector field $\xi$
    is called  {\it  minimal} if $ \xi(M)$ is a minimal submanifold with
    respect to the induced  metric  \cite{G-Z, GM} and  {\it totally
    geodesic} if $\xi(M)$ is a totally geodesic submanifold in $T_1M$. A number of
    examples of locally minimal unit vector fields has been produced by
    J.C.~Gonz\'alez-D\'avila  and L.~Vanhecke \cite{GD-V1}. Most of their examples
    belong to a class of unit vector fields with a non-integrable orthogonal distribution $\xi^\perp$ (the so-called {\it non-holonomic}
    vector fields). The {\it holonomic} case has
    been treated by E.~Boeckx and L.~Vanhecke \cite{BxV1, BxV2} and new examples of
    minimal and harmonic unit vector fields have been produced.

    The {\it totally geodesic } property of the vector field is much more restrictive and
    allows to give a complete description of the field and  the supporting manifold at least for
    2-dimensional manifolds of constant curvature \cite{Ym3}.

    In this paper we treat the case of holonomic unit vector fields,
    namely, the field of unit normals of a given a Riemannian transversally oriented (local)
    hyperfoliation. The question is, \emph{What is the
    connection between the extrinsic geometry of the leaves and the extrinsic (intrinsic) geometry
    of the submanifold $\xi(M)\in T_1M$?} In this paper we give, to some extent, the answer
    to this posed question.

    \section{The results}

    Let $M^{n+1}$ be a Riemannian manifold admitting a transversally oriented  Riemannian
    (local) hyperfoliation. This means that there exists a unit vector field $\xi$ on
    $M^{n+1}$ such that
    the distribution $\xi^\perp$ is integrable, the leaves of the hyperfoliation (~the integral
    submanifolds of $\xi^\perp$ ) are equidistant and the integral trajectories
    of the field $\xi$ are geodesics of $M^{n+1}$. The principal technical  result
    is contained in {\bf Lemma \ref{LeafM}}, which gives an expression for the second fundamental
    form of $\xi(M)$ in terms of the second fundamental forms of leaves and the curvature
    tensor of $M$.  As an application of Lemma \ref{LeafM} to the case of totally umbilic
    hyperfoliation, we have the following
    \vspace{1ex}

    {\bf Theorem \ref{main}.} \ {\it
    Let $\xi$ be a unit normal  vector field of Riemannian transversally orientable
    totally umbilical (local) hyperfoliation on a Riemannian manifold $M$. Then $\xi(M)$ is totally geodesic in $T_1M$ if and only if
    $$
    K_\sigma=\frac{2\,k^2}{k^2-1},
    $$
    where $k=k(s)$ is the value of umbilicity of a leaf and the $K_\sigma$
    are the eigenvalues of the normal Jacobi operator $R(\cdot,\xi)\xi$.
    }
    \vspace{1ex}

    The non-flat spaces of constant curvature evidently drop out from our
    considerations since the $K_\sigma$ are constant along each geodesic. A similar
    property is inherent to all locally symmetric spaces.
    So, the curvature of the manifold should help the vector field to be
    totally geodesic and the manifold has to be non-symmetric. Manifolds and vector fields with
    these desirable properties exist. We completely describe
    2-dimensional Riemannian manifolds admitting {\it geodesic} unit vector fields
    with {\it totally geodesic} property.
    \vspace{1ex}

    {\bf Theorem \ref{2-dim}.} \ {\it
    Let $\xi$ be a (local) unit geodesic vector field on a 2-dimensional Riemannian manifold
    $M$. Then $\xi$ is a totally geodesic vector
    field if and only if the local expression for the metric of $M$ with respect to
    a {\rm(\,$\xi,\xi^\perp$)}-orthogonal coordinate system takes the form
    $$
    ds^2=\frac{(t^2-1)^2}{t^4(t^2+1)^2}\,dt^2+\frac{a^2t^2}{(t^2+1)^2}\,dv^2,
    $$
    where $t$ is the geodesic curvature of $\xi^\perp$-curves, $\xi$ is the
    normalized vector field $\partial_t$ and $a$ is a parameter.
    }
    \vspace{1ex}

    Moreover, we produce an {\it explicit example of a surface of revolution}
    carrying that kind of metric. Let us remark, that the curvature of this surface is
    non-constant, positive for $t^2>1$ and negative for $0<t^2<1$. We also give
    the multidimensional generalization of this example.

\section{Preliminaries.}

    Let $(M,g)$ be an $(n+1)$-dimensional Riemannian manifold  with metric $g$. Denote by
    $\big<\cdot,\cdot\big>$ a scalar product with respect to $g$ and by $\nabla$ the
    Levi-Civita connection on $M$.
    The {\it Sasaki metric} on $TM$ is defined by the following scalar product:
    if $\tilde X,\tilde Y \in TTM$,  then
\begin{equation}
\label{lab5}
        \big<\big< \tilde X,\tilde Y \big>\big>=
        \big<\pi_* \tilde X, \pi_* \tilde Y\big>+\big<K \tilde X,K \tilde Y\big>
\end{equation}
    where $\pi_*:TTM \to TM $ is the differential of the projection $\pi:TM \to M $ and
   $K: TTM \to TM$  is the {\it connection map}.

     Let $\xi$ be a unit vector field on $M$. A vector field $\tilde X \in TTM$ is
     tangent to $\xi(M)$ if and only if \cite{Ym1}
     $$
     \tilde X=(\pi_*\tilde X)^h + (\nabla_{\pi_*\tilde X}\xi)^v,
     $$
    where $(\cdot)^h$ and $(\cdot)^v$ mean {\it horizontal} and {\it vertical} lifts
    of fields into the tangent bundle.

    Introduce a {\it shape operator} $A_\xi$ for the field $\xi$ by
    $$
    A_\xi X=-\nabla_X\xi,
    $$
    where $X$ is an arbitrary vector field on $M$.
    Define a {\it conjugate shape operator} $A_\xi^*$ by
    \begin{equation}\label{conj}
    \big<A_\xi^*X,Y\big>=\big<X,A_\xi Y\big>.
    \end{equation}

    Applying standard singular decomposition of the operator (matrix) $A_\xi$, one can find
    orthonormal local frames $e_0, e_1, \dots , e_n $ and
    $f_0=\xi,\, f_1, \dots , f_n $ on $M$ such that
    $$
        A_\xi \, e_0=0 , \ A_\xi\, e_\alpha =\lambda_\alpha f_\alpha,
        \quad A_\xi^*\, f_0=0,\ A_\xi^*\,f_\alpha=\lambda_\alpha e_\alpha ,
        \ \alpha=1, \dots , n ,
    $$
    where $\lambda_1\geq \lambda_{2}\geq \dots  \lambda_n\geq0$ are real-valued
    functions.

    The frames $e_0, e_1, \dots , e_n $ and $f_0=\xi,\, f_1, \dots , f_n $ are called
    \emph{left} and \emph{right} singular frames respectively for the operator $A_\xi$.

    Remark that one may use if necessary  the signed singular values fixing the directions of the vectors of the singular frame. Setting $\lambda_0=0$, we may rewrite the relations
    on singular frames in a unified form
    \begin{equation}\label{sing}
    \begin{array}{l}
     \ds   A_\xi\, e_i =\lambda_i f_i,\quad
        A_\xi^*\,f_i=\lambda_i e_i ,
        \qquad i=0,1, \dots , n ,\\[2ex]
     \ds   \lambda_0=0, \ \lambda_1,\dots,\lambda_n\geq0.
    \end{array}
    \end{equation}

    The following lemma is easy to prove using (\ref{conj}) and (\ref{sing}).
\begin{lemma}\cite{Ym1}
    At each point of $\xi(M)\subset TM$ the orthonormal frames
    \begin{equation}\label{tang}
    \begin{array}{ll}
     \ds\tilde e_i = \frac{1}{\sqrt{1+\lambda_i^2}}(e_i^h - \lambda_i f_i^v),
        &\ds i=0,1,\dots , n , \\[3ex]
    \ds \tilde n_{\sigma |} =\frac{1}{\sqrt{1+\lambda_\sigma^2}}\big(\lambda_\sigma
    e_\sigma^h +f_\sigma^v \ \big),
    &\ds \sigma=1,\dots , n
    \end{array}
    \end{equation}
    form orthonormal frames in the tangent space of $\xi(M)$ and in the normal space
    of $\xi(M)$, respectively.
\end{lemma}

    Introduce a {\it half tensor } of Riemannian curvature as
    \begin{equation}\label{half}
    r(X,Y)\xi=\nabla_X\nabla_Y\xi-\nabla_{\nabla_XY}\xi.
    \end{equation}

    Now we are able to formulate a lemma, basic for our considerations.

\begin{lemma}\cite{Ym1}\label{Form}
    The components of the second fundamental form of
    $\xi(M)\subset T_1M$ with respect to the frame (\ref{tang}) are given by
    $$
\begin{array}{ll}
    \tilde \Omega_{\sigma | i j}= &\frac{1}{2}\Lambda_{\sigma i j}
    \Big\{\big< r(e_i,e_j)\xi+  r(e_j,e_i)\xi, f_\sigma \big>+\\[2ex]
    &\hspace{3cm}\lambda_\sigma\left[ \lambda_j \big< R(e_\sigma, e_i)
    \xi, f_j \big> +  \lambda_i  \big<R(e_\sigma, e_j) \xi, f_i \big>\right] \Big\},
\end{array}
   $$
    where  $\Lambda_{\sigma i j}=[(1+\lambda_\sigma^2)(1+\lambda_i^2)(1+\lambda_j^2)]^{-1/2}$
    \ $(i,j=0,1,\dots, n;\,\sigma=1,\dots,n)$.
\end{lemma}

\section{The case of Riemannian hyperfoliation}

    Let $M^{n+1}$ be a Riemannian manifold admitting a transversally oriented  Riemannian
    hyperfoliation. In this case $A_\xi X=-\nabla_X \xi$
    is a self-adjoint linear operator on $\xi^\perp$ and for all $X\in\xi^\perp$ it is a shape
    operator of the corresponding leaf. Remark that $A_\xi\,\xi=0$.

    Denote by $\nabla^F$ the induced connection on each leaf. Denote by $B_\xi(X,Y)$
    the second fundamental forms of the leaf, i.e.
    $$
    B_\xi(X,Y)=\big<A_\xi X,Y\big>_F
    $$
    where $\big<\cdot,\cdot\big>_F$ means scalar product with respect to the induced
    metric on each leaf and $X,Y$ are tangent to the corresponding leaf.

    Let $e_\alpha$ \ ($\alpha=1,\dots,n$) be an orthonormal frame consisting of
    eigenvectors of the shape operator, i.e.
    $$
    A_\xi\,e_\alpha=k_\alpha\,e_\alpha,
    $$
    where $k_\alpha$ are the principal curvatures of the corresponding leaf. Since
    $A_\xi$ is self-adjoint, in our notations we have
    $$
    f_0=e_0=\xi, \ f_\alpha=\mathop{\rm sign}(k_\alpha) e_\alpha, \ \lambda_\alpha=|k_\alpha| \quad
    (\alpha=1,\dots,n),
    $$
    where $f_\alpha$ are the vectors of the left singular frame and $\lambda_\alpha$
    are the corresponding singular values. We may simplify notations, if we  set
    \begin{equation}\label{new}
    f_0=e_0=\xi, \ f_\alpha=e_\alpha, \ \lambda_\alpha=k_\alpha \quad
    (\alpha=1,\dots,n),
    \end{equation}
    letting $\lambda_\alpha$ to be not necessarily positive.
    So, the framing of $\xi(M)$ for the
    case under consideration obtains the form
    \begin{equation}\label{fr}
    \begin{array}{ll}
    \tilde e_0=\xi^h,&
    \tilde e_\alpha=\frac{1}{\sqrt{1+k_\alpha^2 }}(e_\alpha^h-k_\alpha\,e_\alpha^v), \\[2ex]
    &\tilde n_{\alpha|}=\frac{1}{\sqrt{1+k_\alpha^2}}(k_\alpha e_\alpha^h+e_\alpha^v).
    \end{array}
    \end{equation}

    Now, the following simplification can be done.

    \begin{lemma}\label{Leaf}
    Let $\xi$ be a unit vector field of Riemannian transversally orientable hyperfoliation
    on a given Riemannian manifold $M^{n+1}$. Denote by $X,Y,Z$ the vector fields
    tangent to the leaf. Then
    $$
    \begin{array}{l}
    \big<r(X,Y)\xi,Z\big>=\big<r(X,Z)\xi,Y\big>=-(\nabla^F_X B_\xi)(Y,Z),\\[1ex]
    \big<r(X,\xi)\xi,Z\big>=-\big<A_\xi X,A_\xi Z\big>_F \\[1ex]
    \big<r(\xi,X)\xi,Z\big>=-\big<A_\xi X,A_\xi Z\big>_F-\big<R(X,\xi)\xi,Z\big>.
    \end{array}
    $$
    \end{lemma}

    \begin{proof}
    Indeed, standard computation yields
    $$
    \big<(\nabla^F_X A_\xi) Y,Z\big>_F= (\nabla^F_X B_\xi)(Y,Z).
    $$
    Consider now $\big<r(X,Y)\xi,Z\big>$. Keeping in mind that $\xi$ is a {\it geodesic}
    vector field, we have
    $$
    \begin{array}{l}
    \big<r(X,Y)\xi,Z\big>=\\[2ex]
    \big<\nabla_X\nabla_Y\xi-\nabla_{\nabla_XY}\xi,Z\big>=
    \big<-\nabla_X(A_\xi Y)-\nabla_{\nabla^F_XY+B_\xi(X,Y)}\xi,Z\big>=\\[2ex]
    \big<-\nabla^F_X(A_\xi Y)-B_\xi(A_\xi Y,X)-\nabla_{\nabla^F_XY}\xi,Z\big>=\\[2ex]
    -\big<(\nabla^F_X A_\xi) Y,Z\big>-\big<A_\xi\nabla^F_X Y,Z\big>+
    \big<A_\xi \nabla^F_XY,Z\big>=-\big<(\nabla^F_X A_\xi) Y,Z\big>.
    \end{array}
    $$
    Thus
    $$
    \big<r(X,Y)\xi,Z\big>=-(\nabla^F_X B_\xi)(Y,Z)=-(\nabla^F_X B_\xi)(Z,Y)=\big<r(X,Z)\xi,Y\big>.
    $$

    Consider $\big<r(X,\xi)\xi,Z\big>$.
    We have
    $$
    \begin{array}{l}
    \big<r(X,\xi)\xi,Z\big>=\\[2ex]
    \big<\nabla_X\nabla_\xi\,\xi-\nabla_{\nabla_X\xi}\xi,Z\big>=
    -\big<\nabla_{-A_\xi X}\,\xi,Z\big>=
    -\big<A_\xi^2X,Z\big>_F=-\big<A_\xi X,A_\xi Z\big>_F.
    \end{array}
    $$
    Finally,
    $$
    \big<r(\xi,X)\xi,Z\big>=\big<r(X,\xi)\xi,Z\big>-\big<R(X,\xi)\xi,Z\big>,
    $$
    which completes the proof.
    \end{proof}

    The following Lemma gives useful information on the relation between extrinsic
    geometry of the leaves of hyperfoliation and extrinsic geometry of the submanifold
    $\xi(M)$ and is a principal tool for further study.

\begin{lemma}\label{LeafM}
Let $\xi$ be a unit normal vector field of Riemannian transversally orientable (local)
hyperfoliation on a given Riemannian manifold $M^{n+1}$. The components of the second
fundamental form of the submanifold $\xi(M)\in T_1M$ with respect to some orthonormal
frame are given by
    $$
    \begin{array}{ll}
    \ds\tilde \Omega_{\sigma|\,00}=&0 \\[2ex]
    \ds\tilde \Omega_{\sigma|\,\alpha
    0}=
    &\frac12\Lambda_{\sigma\alpha 0}
    \Big\{\big[(k_\sigma^2-1)K_{\sigma}-2k_\sigma^2\big]\delta_{\sigma\alpha}
    - \\[1ex]
    &\hspace{4.5cm}\ds(1-k_\alpha k_\sigma)(1-\delta_{\sigma\alpha})\big<R(e_\alpha,\xi)\xi,e_\sigma\big>
    \Big\} \\[2ex]
    \ds\tilde \Omega_{\sigma|\,\alpha\beta}=
    &\frac12\Lambda_{\sigma\alpha\beta}
    \Big\{-2\,(\nabla^F_{e_\sigma}B_\xi)(e_\alpha,e_\beta)
    + \\[1ex]
    &\hspace{1cm}\ds(1-k_\sigma k_\alpha)\big<R(\xi,e_\alpha)e_\beta,e_\sigma\big>+
    \ds(1-k_\sigma k_\beta)\big<R(\xi,e_\beta)e_\alpha,e_\sigma\big>
    \Big\},
    \end{array}
    $$
    where $K_\sigma$ are the eigenvalues of the normal Jacobi operator
    $R(\cdot,\xi)\xi$ and $\delta_{\sigma\alpha}$ is the Kronecker symbol.
    \end{lemma}

    Note that Lemma \ref{LeafM} can be applied to the case of a local foliation such as a
    family of distance spheres, tubes etc. As an immediate corollary we see that if the leaves
    are totally geodesic or even totally umbilic then $\xi(M)$  is a \emph{minimal} submanifold
    but it is not totally geodesic in general.

\begin{proof}
    (a)\ Since $\xi$ is a geodesic vector field, we may set $e_0=\xi$ and therefore we have
    $\big<r(e_0,e_0)\xi,e_\sigma\big>=0$.
    Applying Lemma \ref{Form}, we get $\tilde \Omega_{\sigma|\,00}=0.$

    (b)\ From Lemma \ref{Form}
    $$
    \tilde\Omega_{\sigma|\,\alpha0}=\frac12\Lambda_{\sigma\alpha 0}
    \Big\{\big<r(e_\alpha, e_0) \xi,f_\sigma \big>+ \big<r(e_0, e_\alpha) \xi, f_\sigma \big>+
    \lambda_\sigma\lambda_\alpha \big< R(e_\sigma, e_0)\xi, f_\alpha \big>\Big\}.
    $$
    Taking into account (\ref{new}) and applying Lemma \ref{Leaf}, we get
    $$
    \begin{array}{l}
    \big<r(e_\alpha, e_0) \xi,f_\sigma \big>=\big<r(e_\alpha, \xi) \xi,e_\sigma
    \big>=-k_\alpha
    k_\sigma\big<e_\alpha,e_\sigma\big>=-k_\sigma^2\delta_{\sigma\alpha}\\[1ex]
    \big<r(e_0,e_\alpha) \xi,f_\sigma \big>=\big<r(\xi,e_\alpha) \xi,e_\sigma \big>=
    -k_\sigma^2\delta_{\sigma\alpha}- \big<R(e_\sigma,\xi)\xi,e_\alpha\big>.
    \end{array}
    $$
    On the other hand, setting $K_\alpha=\big<R(e_\alpha,\xi)\xi,e_\alpha\big>$, we
    have
    $$
    \big<R(e_\sigma,\xi)\xi,e_\alpha\big>=K_{\sigma}\delta_{\sigma\alpha}+
    (1-\delta_{\sigma\alpha})\big<R(e_\sigma,\xi)\xi,e_\alpha\big>.
    $$
    After substitutions, we get
    $$
    \ds\tilde \Omega_{\sigma|\,\alpha
    0}=\frac12\Lambda_{\sigma\alpha 0}
    \Big\{\big[(k_\sigma^2-1)K_{\sigma}-2k_\sigma^2\big]\delta_{\sigma\alpha}
    - \ds(1-k_\alpha k_\sigma)(1-\delta_{\sigma\alpha})\big<R(e_\alpha,\xi)\xi,e_\sigma\big>
    \Big\}
    $$

    (c) \ From Lemma \ref{Form} and (\ref{new})
    $$
    \begin{array}{rcl}
    \tilde \Omega_{\sigma |\,\alpha\beta}
    &=& \frac{1}{2}\Lambda_{\sigma\alpha\beta}
    \Big\{ \big<
    r(e_\alpha, e_\beta) \xi, e_\sigma \big>+ \big<r(e_\beta, e_\alpha) \xi, e_\sigma \big>\\[2ex]
    &&+ k_\sigma\left[ k_\alpha \big< R(e_\sigma, e_\beta)
    \xi, e_\alpha \big> +  k_\beta  \big<R(e_\sigma, e_\alpha) \xi, e_\beta \big>\right] \Big\},
    \end{array}
    $$
    Lemma \ref{Leaf} and the Codazzi equation yield
    $$
    \begin{array}{l}
    r(e_\alpha, e_\beta) \xi, e_\sigma \big>=-(\nabla^F_{e_\alpha}
    B_\xi)(e_\beta,e_\sigma),\quad
    r(e_\beta, e_\alpha) \xi, e_\sigma \big>=-(\nabla^F_{e_\beta}
    B_\xi)(e_\alpha,e_\sigma),\\[2ex]
    \big< R(e_\sigma, e_\beta)\xi, e_\alpha \big>=-(\nabla^F_{e_\sigma}
    B_\xi)(e_\beta,e_\alpha)-\big<r(e_\beta,e_\alpha)\xi,e_\sigma\big>,\\[2ex]
    \big< R(e_\sigma, e_\alpha)\xi, e_\beta \big>=-(\nabla^F_{e_\sigma}B_\xi)(e_\alpha,e_\beta)
    -\big<r(e_\alpha,e_\beta)\xi,e_\sigma\big>.
    \end{array}
    $$
    So we have
    $$
    \begin{array}{l}
    \big<r(e_\alpha, e_\beta) \xi, e_\sigma \big>+
    \big<r(e_\beta, e_\alpha) \xi, e_\sigma \big>=\\[2ex]
    \qquad \qquad -2(\nabla^F_{e_\sigma}B_\xi)(e_\alpha,e_\beta)-
    \big< R(e_\sigma, e_\alpha)\xi, e_\beta \big>-
    \big< R(e_\sigma, e_\beta)\xi, e_\alpha \big>.
    \end{array}
    $$
    After substitutions, we get
    $$
    \begin{array}{ll}
    \ds\tilde \Omega_{\sigma|\,\alpha\beta}=
    &\frac12\Lambda_{\sigma\alpha\beta}
    \Big\{-2(\nabla^F_{e_\sigma}B_\xi)(e_\alpha,e_\beta)
    + \\[1ex]
    &\hspace{1cm}\ds(1-k_\sigma k_\alpha)\big<R(\xi,e_\alpha)e_\beta,e_\sigma\big>+
    \ds(1-k_\sigma k_\beta)\big<R(\xi,e_\beta)e_\alpha,e_\sigma\big>
    \Big\},
    \end{array}
    $$
    which completes the proof.
\end{proof}

    Now we can characterize the totally umbilic foliations as follows.
\begin{theorem} \label{main}
    Let $\xi$ be a unit normal  vector field of Riemannian transversally orientable
    totally umbilical (local) hyperfoliation on a Riemannian manifold $M$. Then $\xi(M)$ is totally geodesic in $T_1M$ if and only if
    \begin{equation}\label{cond}
    K_\sigma=\frac{2\,k^2}{k^2-1},
    \end{equation}
    where $k=k(s)$ is the value of umbilicity of a leaf and $K_\sigma$
    are the eigenvalues of the normal Jacobi operator $R(\cdot,\xi)\xi$.
    \end{theorem}

\begin{proof}

    The result of Lemma \ref{LeafM} means that the extrinsic geometry of holonomic, i.e.
    with integrable distribution $\xi^\perp$,
    geo\-de\-sic vector fields depends on the extrinsic geometry of leaves and the normal
    Jacobi operator $R_\xi=R(\cdot,\xi)\xi$. A submanifold $F\subset M$ is said to be {\it
    curvature adapted} \cite{BV} if for every normal vector $\xi$ to $F$ at a point $p\in
    F$ the
    following conditions hold:
    $$
    \begin{array}{l}
     R_\xi(T_p F)\subset T_pF \\[1ex]
     A_\xi\circ R_\xi=R_\xi\circ A_\xi,
    \end{array}
    $$
    where $A_\xi$ is the shape operator of $F$. The first condition is always fulfilled
    for a hypersurface. The second means that there exists a basis of $T_pF$ consisting
    of eigenvectors of both  $R_\xi$ and $A_\xi$. Every totally umbilical submanifold
    is curvature adapted and has parallel second fundamental form. These
    facts immediately imply
    \begin{proposition}\label{totumbil}
    Let $\xi$ be a unit normal  vector field of Riemannian transversally orientable
    totally umbilical hyperfoliation on a given Riemannian manifold $M^{n+1}$.
    The non-zero components of the second fundamental  form
    for $\xi(M)\in T_1M$ are given by
    $$
    \tilde\Omega_{\sigma|\,\sigma0}=\frac12\,\frac{1}{1+k^2}\Big[(k^2-1) K_{\sigma}-2k^2
    \Big],
    $$
    where $k=k(s)$ is the value of umbilicity of a leaf $F^n_s$ and $K_\sigma$ are the
    eigenvalues of the normal Jacobi operator $R(\cdot,\xi)\xi$.
    \end{proposition}

Now the main result follows immediately from Proposition \ref{totumbil}.

    \end{proof}

    Remark that $K_\sigma= K_{\xi\wedge e_\sigma}$, where $K_{\xi\wedge e_\sigma}$ means a sectional curvature along the plane $\xi\wedge e_\sigma$ and
    $e_\sigma$ are the eigenvectors of the normal Jacobi operator. A similar condition
    is necessary for the totally geodesic property in the case of curvature adapted foliation
    (even a local one), namely
        $$
    K_\sigma=\frac{2\,k_\sigma^2}{k_\sigma^2-1}.
    $$
    This condition fails if $M^{n+1}$ is locally symmetric and the leaves are
    homogeneous. In this case $K_\sigma$ are constant along $\xi$-geodesics while
    $k_\sigma$ are the functions of its natural parameter. Typical examples are provided by
    the field of unit normals  of a family of geodesic spheres or by the tubes around a totally
    geodesic submanifold. These vector fields are minimal \cite{BxV1} but never totally geodesic.

    As a direct application of Lemma \ref{LeafM} to the case of a 2-dimensional Riemannian
    manifold, we are able to describe completely  the totally geodesic unit vector fields belonging
    to the class under consideration  and the supporting manifold in the following
    terms.

    \begin{theorem}\label{2-dim}
    Let $\xi$ be a (local) unit geodesic vector field on a 2-dimensional Riemannian manifold
    $M$. Then $\xi$ is a totally geodesic vector
    field if and only if the local expression for the metric of $M$ with respect to a
    {\rm(\,$\xi,\xi^\perp$)}-orthogonal coordinate system takes the form
    \begin{equation}\label{metric}
    ds^2=\frac{(t^2-1)^2}{t^4(t^2+1)^2}\,dt^2+\frac{a^2t^2}{(t^2+1)^2}\,dv^2,
    \end{equation}
    where $t$ is the geodesic curvature of $\xi^\perp$-curves, $\xi$ is the
    normalized vector field $\partial_t$ and $a$ is a parameter.
    \end{theorem}

    \begin{proof}
     Let $\xi$ be a (local) geodesic unit vector field on a 2-dimensional Riemannian
    manifold $M$ of Gaussian curvature $K$. The result of Lemma \ref{LeafM}
    allows to simplify the matrix of the second fundamental form of $\xi(M)\in T_1M$
    to
    $$
    \Omega=\left(
    \begin{array}{cc}
    \ds0&\ds\frac12\frac{(k^2-1)K-2k^2}{1+k^2}\\[3ex]
    \ds\frac12\frac{(k^2-1)K-2k^2}{1+k^2}& -\ds\frac{e_1(k)}{(1+k^2)^{3/2}}
    \end{array}
    \right),
    $$
    where $k$ is the geodesic curvature of the integral trajectories of the unit vector
    field $e_1=\xi^\perp$.

    Taking $\xi$-integral trajectories as the first family of coordinate lines and
    $e_1=\xi^\perp$-integral trajectories as the second one,
    we can express the metric of $M^2$ in the form
    $$
    ds^2=du^2+g^2(u,v)\,dv^2,
    $$
    where $g(u,v)$ is some (positive) function. Remark that the geodesic curvature of
    $e_1$-curves with respect to our coordinate system takes the form
    \begin{equation}\label{geocurv}
    k=-\frac{g_u}{g}.
    \end{equation}
    Suppose now that $\xi(M)$ is totally geodesic in $T_1M$. Then $e_1(k)=0$ and hence
    $k$ does not depend on the $v$-parameter. Solving (\ref{geocurv}) with respect to $g$, we get
    $$
    g(u,v)=C(v)\exp (-\int k(u)\,du).
    $$
    After $v$-parameter change, we reduce the metric to the form of metric of a
    sur\-face of revolution
    $$
    ds^2=du^2+f^2(u)\,dv^2.
    $$
    So the curves $u=const$ (the parallels) give us a totally umbilical
    foliation on $M^2$. The \emph{value of umbilicity} is the geodesic curvature of
    parallels, the vector field $\xi=\partial_u$ is a unit vector field tangent to meridians,
    the Gaussian curvature $K$ of $M^2$ for this depends only on $u$ -- the natural parameter
    on meridians -- and  $K=k'(u)-k^2(u)$. To satisfy the totally geodesic property,
    the geodesic curvature $k$ has to be a solution of the differential equation
    $$
    k'=\frac{k^2\,(k^2+1)}{k^2-1}.
    $$
    The implicit solution is $ \ds u=2\arctan k+\frac{1}{k}+u_0.$
    The inverse function $k=k(u)$ exists on intervals where $k(u)\ne1$.

    To produce an explicit solution,  we proceed as follows.
    Choose $k$ as a parameter, say $t$. Since $\ds k=-f^{-1}f'_u$ , we can write two
    relations
    \begin{equation}\label{func}
    \frac{f'_u}{f}=-t\, ,\quad \quad \frac{d\,t}{du}=\frac{t^2(t^2+1)}{t^2-1}.
    \end{equation}
    Making a parameter change, we obtain a differential equation on $f (t)$ of the form
    $\ds
    \frac{f'_t}{f}=\frac{1-t^2}{t(1+t^2)}
    $
    with a general solution
    $\ds
    f(t)=\frac{a\,t}{t^2+1},
    $
    where $a$ is the constant of integration. From (\ref{func}) we can also find
    $\ds
    du=\frac{t^2-1}{t^2(t^2+1)}\, dt
    $
    and therefore the metric under consideration takes the form (\ref{metric})
    with respect to the parameters $(t,v)$.

    \end{proof}

    Indeed, we are able to get an isometric immersion of the metric constructed
     into Euclidean 3-space as a {\it surface} of revolution.
    Some additional considerations show that we get the most regular surface for $a=1$.
    \vspace{1ex}

    \begin{example}
    Let $\big\{x(t),z(t)\big\}$ be a profile curve, generating a surface with the metric
    (\ref{metric}) with $a=1$. Then, evidently,
    $\ds
    x(t)=\frac{t}{(t^2+1)}
    $
    and
    $\ds
    (x'_t)^2+(z'_t)^2=\frac{(t^2-1)^2}{t^4(t^2+1)^2}.
    $
    From this we find $$\ds z'_t=\pm \frac{(t^2-1)\sqrt{2t^2+1}}{t^2(t^2+1)^2}.$$
    Choose one branch, say with a positive sign. A relatively simple calculation gives
    $\ds z(t)=\frac{(2t^2+1)^{3/2}}{t(t^2+1)}$ (up
    to an additive constant). Thus, finally, we have a parametric curve
    $$
    x(t)=\frac{t}{(t^2+1)}, \quad
    z(t)= \frac{(2\,t^2+1)^{3/2}}{t(t^2+1)}
    $$
    parametrized with the geodesic curvature of the meridians of the associated  surface of
    revolution
    and having one singular point corresponding to $t=1$. The following picture gives a
    graph.

    \begin{picture}(0,0)
    \put(100,0){\special{em:graph grafic.bmp}}
    \end{picture}
    \vspace{3.5cm}

    Remark that $ K<0 $ for $ t^2\in(0,1)$ and $ K>0$  for $t^2\in(1,+\infty)$.
    The point $(0,0,2\sqrt{2})$ is an umbilical one at infinity (for given parameterization).

    \end{example}

    The example is not essentially 2-dimensional. Consider a  metric of revolution of the form
        $$
        ds^2=du^2+f^2(u)\sum_{\alpha=1}^n (d v^\alpha)^2 .
        $$
        Then the leaves of hyperfoliation $u=const$ are all totally umbilic with a value of umbilicity $k=-f^{-1}f'_u$. To make the vector field $\xi=\partial_u$ totally geodesic, this value should satisfy the same differential equation, namely
        $$
    k'=\frac{k^2\,(k^2+1)}{k^2-1},
    $$
which has the same solution as in 2-dimensional example.

\vspace{1cm}

\noindent
Department of Geometry,\\
Faculty of Mechanics and Mathematics,\\
Kharkiv National University,\\
Svobody Sq. 4,\\
 61077, Kharkiv,\\
Ukraine.\\
e-mail: yamp@univer.kharkov.ua

\end{document}